\documentclass{article}
\usepackage{arxiv}

\usepackage{tabularx}
\usepackage{amsmath,amssymb,amsfonts}
\usepackage{textcomp}
\usepackage{accents}
\usepackage{graphicx,color} %
\usepackage{epsfig} %
\usepackage{times} %
\usepackage[T1]{fontenc}
\usepackage[utf8]{inputenc} %
\usepackage{makecell}
\usepackage{algorithm}
\usepackage[noend]{algpseudocode}

\usepackage{enumitem}
\usepackage{mathtools}
\usepackage{multirow}

\usepackage{hyperref}

\graphicspath{{figures/}}

\newcommand{\rr}{\mathbb{R}}

\newcommand{\nn}{\mathbb{Z}}

\newcommand{\eqdef  }{\triangleq}

\newcommand{\II}{{\mathcal{I}}}

\newcommand{\minus}{\scalebox{0.75}[1.0]{$-$}}
\newcommand{\matrice}[2]{\left[\hspace*{-.1cm}\begin{array}{#1} #2 \end{array}\hspace*{-.1cm}\right]}

\newcommand{\FIT}{\mathrm{FIT}}

\newcommand{\atan}{\mathop{\rm atan}\nolimits}
\newcommand{\NRMSE}{\mathop{\rm NRMSE}\nolimits}

\newcommand{\Nsamp}{N}

\renewcommand{\headeright}{}
\renewcommand{\undertitle}{}

\title{A machine-learning approach to synthesize virtual sensors for  parameter-varying systems}

\author{Daniele Masti \thanks{Corresponding author. E-mail address: \texttt{daniele.masti@imtlucca.it}}\\
IMT School for Advanced Studies Lucca, Piazza San Francesco 19, 55100 Lucca, Italy
\And
Daniele Bernardini\\
ODYS S.R.L., Via Pastrengo, 14, 20159 Milano, Italy
\And
Alberto Bemporad\\
IMT School for Advanced Studies Lucca, Piazza San Francesco 19, 55100 Lucca, Italy
}

\begin{document}
\maketitle
\begin{abstract}
This paper introduces a novel model-free approach to synthesize virtual sensors
for the estimation of dynamical quantities that are unmeasurable at runtime but 
are available for design purposes on test benches. After collecting a dataset 
of measurements of such quantities, together with other variables that are 
also available during on-line operations, 
the virtual sensor is obtained using machine learning techniques 
by training a predictor whose inputs are the measured variables and the features extracted by 
a bank of linear observers fed with the same measures. The approach is applicable
to infer the value of quantities such as physical states and other time-varying 
parameters that affect the dynamics of the system. The proposed virtual sensor 
architecture --- whose structure can be related to the Multiple Model Adaptive 
Estimation framework --- is conceived to keep computational and memory 
requirements as low as possible, so that it can be efficiently implemented in 
embedded hardware platforms. 

The effectiveness of the approach is shown in different numerical examples, 
involving the estimation of the scheduling parameter of a nonlinear 
parameter-varying system, the reconstruction of the mode of a switching linear 
system, and the estimation of the state of charge (SoC) of a lithium-ion 
battery.
\end{abstract}


\section{Introduction}

Most real-world processes exhibit complex nonlinear dynamics that are difficult 
to model, not only because of nonlinear interactions between input and output 
variables, but also because of the presence of time-varying signals that change 
the way the involved quantities interact over time. A typical instance is the 
case of systems subject to wear of components, in which the dynamics slowly 
drift from a nominal behavior to an aged one, or systems affected by 
slowly-varying unknown disturbances, such as unmeasured changes of ambient 
conditions. Such systems can be well described using a \emph{parameter-varying} 
model~\cite{lopez2019review} that depends on a vector $\rho_k \in \rr^S$ of 
parameters, that in turn evolves over time:
\begin{equation}
\Sigma_{P}\eqdef \ \left\{
\begin{array}{cl}
x_{k+1}=  f(x_k,u_k,\rho_k)\\
\rho_{k+1}=h(\rho_k,k,u_k)\\
y_k= g(x_k,\rho_k)\\
\end{array}
\right.
\label{eq:true-model}
\end{equation}
where $x_k \in \rr^{n_x}$ is the state vector, $y_k \in \rr^{n_y}$ is the output vector, $u_k \in \rr^{n_u}$ is the input vector, 
$f:\rr^{n_x} \times \rr^{n_u}\times \rr^{S}\to\rr^{n_x}$, $g:\rr^{n_x}\times 
\rr^{S}\to\rr^{n_y}$ and $h:\rr^{S}\times\rr^{n_u}\times\rr \to \rr^{S}$. In 
this paper we assume that the mappings in~\eqref{eq:true-model} are  
\textit{unknown}.

Special cases of~\eqref{eq:true-model} widely studied in the literature are linear parameter-varying (LPV) systems~\cite{toth2010LPVModeling}, in which $f$, $g$ are linear functions of $x_k$, $u_k$, and switched affine systems~\cite{Torrisi2004HYSDEL}, in which $\rho_k$ only assumes a value within a finite set.

Inferring the value of $\rho_k$ in real time from input/output data can be useful for several reasons. In predictive maintenance and anomaly/fault detection~\cite{Rotondo2013FTC,Witczak2014faultDiagnosis,guzman2021actuator}, detecting a drift in the value of $\rho_k$ from its nominal value or range of values can be used first to detect a fault and then to isolate its nature. In gain-scheduling control~\cite{slotine1991applied,misin2020lpv}, $\rho_k$ can be used instead to decide the control law to apply at each given time instant. 

Due to the importance of estimating $\rho_k$, various solutions have been 
proposed in the literature to estimate it during system operations. 
If the model in~\eqref{eq:true-model} were known, even if only approximately, 
nonlinear and robust state estimators could be successfully 
applied~\cite{do2020h}. On the other 
hand, if the mechanism regulating the interaction between $\rho_k$ and the 
measurable quantities (usually $u_k$ and $y_k$) is not known, but a 
dataset of historical data is available, the classical indirect approach would 
be to identify an overall model of $\Sigma_P$ using nonlinear system identification 
techniques~\cite{keesman2011system,masti2018Autoencoder} and then build a 
model-based observer to estimate $\rho_k$. The drawbacks of such an indirect 
approach are that it can be a very time-consuming task and that the resulting 
model-based observer can be complex to implement. This issue is especially cumbersome if one is 
ultimately interested in just getting an observer and have no further use  for the model itself.

\textit{Virtual sensors}~\cite{milanese2009filter,Poggi2012HighSpeedPiecewise} provide an alternative approach to solve such a problem: the idea is to build an \emph{end-to-end} estimator for $\rho_k$ by directly learning from data the mapping from measured inputs and outputs to $\rho_k$ itself.
The approach is interesting because it does not require identifying a full model of the system from data, nor it requires simplifying an existing model (such as a high-fidelity simulation model) that would be otherwise too complex for model-based observer design. Similar estimation problems have been tackled in the context of novelty detection~\cite{marchi2015novel} and of time-series clustering~\cite{aghabozorgi2015time,morse2007efficient}. 

\subsection{Contribution}
\label{sec:contribution}

The goal of this paper is to develop an approach to synthesize virtual 
sensors that can estimate $\rho_k$ when its measurements are not available by using 
data acquired when such a quantity is directly measurable. Such a scenario
often arises in serial production, in which the cost of components must be 
severely reduced. The purpose of the proposed approach is to enable replacing 
physical sensors with lines of code.

The method developed in this paper is loosely related to Multiple 
Model Adaptive Estimation~\cite{akca2019multipleModelSurvey} (MMAE)
and consists of three main steps: 
\begin{enumerate}
\item Learn a finite set of simple linear time-invariant (LTI) models from data 
that roughly covers the behavior of the system for the entire range of values 
of $\rho_k$ of interest;
\item Design a set of standard linear observers based on such models;
\item Use machine-learning methods to train a lightweight predictor that maps 
the estimates obtained by the 
observers and raw input/output signals into an estimate $\hat \rho_k$ of 
$\rho_k$.
\end{enumerate}
To do so, this paper extends the preliminary results presented 
in~\cite{Masti2019VirtualSensor} in several ways: it formulates the problem for
nonlinear systems; it explores the performance of the approach for 
mode-discrimination of switching systems and it provides a thorough performance 
analysis of various lightweight machine-learning techniques that can be used to 
parameterize the virtual sensor architecture. In doing so, it also provides an entirely off-line alternative strategy for identifying the local linear models required to synthesize the bank of observers based on an interpretation of well-known decision tree regressors as a supervised clustering scheme.

The intuition behind our approach is that, in many cases of practical interest, 
the dynamics of $\rho_k$ are slower than the other dynamics of the system. 
This fact suggests that a linear model identified on a dataset in which $\rho_k$ is 
close to a certain value $\bar \rho$ will well approximate $\Sigma_P$ for all 
$\rho_k\approx\bar \rho$. Following this idea, we envision a scheme in which $N_\theta$ values $\bar\rho_i$, 
$i=1,\dots, N_\theta$ are automatically selected and, for each value $\bar \rho_i$, 
a linear model is identified and a corresponding linear observer synthesized.
A machine-learning algorithm is then used to train a predictor that consumes 
the performance indicators constructed from such observers, together with raw 
input and output data, to produce an estimate $\hat \rho_k$ of $\rho_k$ at each 
given time $k$.

The paper is organized as follows: in Section~\ref{sec:MMAE} we recall the MMAE 
framework and introduce the necessary steps to bridge such a model-based 
technique to a data-driven framework. In Section~\ref{sec:virtualSensorArch}, 
we detail the overall virtual sensor architecture and the internal structure 
of its components. Section~\ref{sec:results} is devoted to studying the quality of 
estimations and the numerical complexity of the synthesized virtual sensor on 
some selected nonlinear and piecewise affine (PWA) benchmark problems, 
including the problem of estimating the state of charge of a battery,  to 
establish both the estimation performance of the approach and the 
influence of its hyper-parameters. Finally, some conclusions 
are drawn in Section~\ref{sec:conclusions}.

The Python code to reproduce the results described in the paper is available at \url{http://dysco.imtlucca.it/masti/mlvs/rep_package.zip}.

\section{Multiple Model Adaptive Estimation}
\label{sec:MMAE}
Following the formulation 
in~\cite{bar2004estimation,alsuwaidan2011generalized}, this section recalls the 
main concept of the MMAE approach. Consider the dynamical system
\begin{equation}
\Sigma \eqdef \  \left\{ \begin{array}{rcl}
x_{k+1}&=&f(\theta_k,x_k,u_k)\\
y_{k}&=&h(\theta_k,x_k,u_k)
\end{array}
\right.
\label{eq:MMAEbaseSystem}
\end{equation}
in which $\theta_k \in \rr^{n_\theta}$ is a generic  parameter vector.
The overall idea of MMAE is to use a bank of $N_\theta$ state 
estimators\footnote{Usually Kalman filters (KF) are employed, but exceptions
exist~\cite{akca2019multipleModelSurvey} .} --- each one associated to a
specific value $\theta_i \in \Theta \coloneqq 
\{\theta_1,\dots,\theta_{N_\theta}\}$ 
--- together with a hypothesis testing algorithm to infer information 
about~\eqref{eq:MMAEbaseSystem}, e.g.: to build an estimate $\hat x_k$ of 
$x_k$. 
In this scheme, the intended purpose of the latter component is to infer, from 
the behavior of each observer, which one among the different models 
(``hypotheses'') is closest to the underlying process, and use such 
information to construct an estimate $\hat x$ of the state of $\Sigma$. For 
linear time-invariant (LTI) 
representations, a classical approach to do so is to formulate the hypothesis 
tester as an appropriate statistical test, exploiting the fact that the 
residual signal produced by a properly matched KF is a zero-mean white-noise 
signal.

MMAE is a model-based technique in that it requires a model of the process, a 
set $\Theta$ of parameter vectors, and a proper characterization of the noise 
signals supposed to act on the system.
Among those requirements, determining $\Theta$ is especially crucial to get reliable results as, at each time, at least one value $\theta_j \in \Theta$ must describe the dynamics of the underlying system accurately enough. In many practical situations, it is not easy to find a good tradeoff between keeping $N_\theta$ large enough to cover the entire range of the dynamics and, at the same time, small enough to limit the computational burden manageable and avoid the tendency of MMAE to work poorly if too many models are considered~\cite{li1996MMAEwithVariableStructure}.
Another difficulty associated with MMAE schemes is the reliance on models to 
synthesize the hypothesis tester. Moreover, many approaches require 
sophisticated statistical arguments, which can hardly be tailored to 
user-specific needs.

\section{Data-driven determination of linear models}
\label{sec:virtualSensorArch}
The first step to derive the proposed data-driven virtual-sensor is to 
reconcile the MMAE framework with the parameter-varying model description 
in~\eqref{eq:true-model}. Assume for the moment that $f$ and $g$ 
in~\eqref{eq:true-model} are known and differentiable. Then, in the 
neighborhood of an arbitrary tuple $(\bar \rho,\bar x,\bar u)$ 
it is possible to approximate~\eqref{eq:true-model} by
\begin{equation}
\begin{array}{rll} x_{k+1}-\bar x& \approx 
f(\bar x,\bar u,\bar \rho)-\bar x+\nabla_x f(\bar x,\bar u,\bar \rho) (x_k-\bar x)+\nabla_u f(\bar x,\bar u,\bar \rho) (u_k-\bar u)\\
y_k& \approx g(\bar x,\bar \rho)+\nabla_x g(\bar x,\bar u,\bar \rho) (x_k-\bar 
x)
\end{array}
\label{eq:linearized-model}
\end{equation}
In~\eqref{eq:linearized-model} the contributions of the
Jacobians with respect to $\rho$ is neglected due to the fact that, as 
mentioned earlier, $\rho_k$ is assumed to move slowly enough to remain close 
to $\bar \rho$ within a certain time interval, meaning the neglected Jacobians would be 
multiplied by $\rho_k-\bar\rho\approx 0$. 
Hence, from~\eqref{eq:linearized-model} we can derive the following affine 
parameter-varying (APV)
approximation of~\eqref{eq:true-model}
\begin{equation}
\begin{array}{rcl}
 x_{k+1}&\approx & A(\rho_k) x_k + B(\rho_k) u_k +d(\rho_k)\\
 y_k &\approx & C(\rho_k) x_k+e(\rho_k)
\end{array}
\label{eq:APV}
\end{equation}
in which the contribution of the constant terms $\bar x$, $\bar u$ is
contained in the bias terms $d(\rho_k)$, $e(\rho_k)$. 
In conclusion, if $\Sigma_P$ were known, a MMAE scheme could be used to compute 
the likelihood that the process is operating 
around a tuple $(\bar x, \bar u, \rho_i)$, where $\rho_i \in 
\Theta^\rho\eqdef  \{\rho_1, \dots,\rho_{N_\theta} \}$ is used
in place of the parameter vector $\theta_k$ in~\eqref{eq:MMAEbaseSystem}.

\subsection{Learning the local models}
\label{sec:linear_models}
As model~\eqref{eq:true-model} is not available, we need to identify
the set of linear (affine) models in~\eqref{eq:APV} from data.
Assuming that direct measurements of the state $x_k$ of the physical system
are not available, we restrict affine autoregressive models with exogenous 
inputs (ARX) of a fixed order, each of them uniquely identified by a parameter 
vector $\gamma\in\rr^{n_\gamma}$.

Learning an APV approximation of $\Sigma_P$ amounts to train a functional 
approximator $M_{LPV}:\rr^{S}\to\rr^{n_\gamma}$ to predict the correct 
vector $\gamma_i$ corresponding
to any given $\bar\rho_i$.  Given a dataset $D_\Nsamp \coloneqq 
\{u_k,y_k,\rho_k\}$, $\ k=1,\dots,\Nsamp$, of samples acquired via an 
experiment on the real process, such a training problem is solved by
the following optimization problem
\begin{equation}
\begin{array}{rrl}
 \underset{M_{LPV}}{\mathrm{min}} & \displaystyle \sum_{k=k_1}^{\Nsamp}& 
 L_{M_{LPV}}(\hat y_{k}, y_{k}) \\
  &\mathrm{subject~to} &\hat y_{k}=[\minus y_{k-M},\dots,\minus y_{k-1},u_{k-M},\dots,u_{k-1}~1]\gamma_{k}\\
  &&\gamma_{k}=M_{LPV}(\rho_k)\\
  &&k=k_1,\dots,\Nsamp
\end{array}
\label{eq:MLPV_learning}
\end{equation}
where $k_1\eqdef  M+1$, and $L_{M_{LPV}}$ is an appropriate loss function. Note 
that, as commonly expected when synthesizing virtual sensors, we assume that 
measurements of $\rho_k$ are available for training, although they will 
not be during the operation of the virtual sensor. Moreover, note that problem 
$M_{LPV}$ is solved offline, so the computation requirements of the regression 
techniques used to solve~\eqref{eq:MLPV_learning} are not of concern.

Compared to adopting a recursive system identification technique to learn a local linear model of the process at each time $k$, and then associate each $\gamma_k$ to its $\rho_k$ (e.g., by using Kalman filtering techniques~\cite{Masti2019VirtualSensor,ljung1987system}), the approach in~\eqref{eq:MLPV_learning} does not require tuning the recursive identification algorithm and takes into account the value of $\rho_k$ at each $k$. This prevents that similar values of $\rho$ are associated with very different values of $\gamma$, assuming that the resulting function $M_{LPV}$ is smooth enough.

\subsubsection{An end-to-end approach to select the representative models}

\label{sec:DTRforSegmentation}

By directly solving~\eqref{eq:MLPV_learning}, a set $\Gamma\eqdef   
\{\gamma_i\}_{i=M+1,\dots,\Nsamp}$ of local models is obtained. Using 
$\Theta=\Gamma^\rho$ in an MMAE-like scheme would 
result in an excessively complex scheme. To address this issue, a smaller set of models could be extracted by running a clustering algorithm on the dataset $\Gamma$, and the set $\Theta^\rho$ of representative models selected as the set of the centroids of the found clusters. Doing so, however%
A better idea comes from observing that some regression techniques, such as decision-tree regressor~\cite{hastie2009elements} (DTRs), naturally produce piece-wise constant predictions, which suggests the following alternative method: ($i$) train a DTR to learn an predictor $\hat M_{LPV}:\rr^{S}\to \rr^{n_\gamma} \times \rr^{S}$ (in an autoencoder-like fashion~\cite{HintonSalakhutdinov2006b,feng2018autoencoderTree}), possibly imposing a limit on its maximum depth; ($ii$) set $\Theta$ as the leaves $\bar\gamma_{j}$ of the grown tree $\hat M_{LPV}$. %

Compared to using a clustering approach like 
K-means~\cite{lloyd1982kmeans,bemporad2018fittingJumpModels}, the use of DTRs 
does not require selecting a fixed number of clusters a priori and also 
actively takes into account the relation between $\rho$ and $\gamma$. 
In fact, with the proposed DTR-based approach, the regression tree will not grow in regions where $\rho$ is not informative about $\gamma$, therefore aggregating a possibly large set of
values of $\gamma$ with the same representative leaf-value $\bar\gamma_j$. This 
latter aspect is important for our ultimate goal of exploiting the resulting 
set of models to build a bank of linear observers.

Once the set $\Theta^\rho=\{\bar\gamma_j\}_{j=1}^{N_\theta}$ of local ARX 
models has been selected, each of them is converted into a corresponding 
minimal state-space representation in observer canonical form~\cite{mellodge2015practical}
\begin{equation}
    \Sigma_j:=\ \left\{\begin{array}{rcl}
        \xi^j_{k+1}&=&A_j\xi^j_k+B_ju_k+d_j\\
        y_k&=&C_j\xi^j_k+e_j
    \end{array}\right.\quad j=1,\ldots,N_\theta
    \label{eq:Sigmaj}
\end{equation}
As all vectors $\bar\gamma_j$ have the same dimension $n^\gamma$, we assume 
that all states $\xi^j$ have the same dimension $v$. The models $\Sigma_j$
in~\eqref{eq:Sigmaj} are used to design a corresponding linear observer, as 
described next.

\subsection{Design of the observer bank}

For each model $\Sigma_j$, we want to design an observer providing an estimate 
$\hat \xi_k^j$ of the state $\xi_k^j$ of $\Sigma_j$. Let $i^j_k \in R^{v}$ be 
the \emph{information vector} generated by 
the observer at time $k$, which includes $\hat \xi_k^j$ and possibly
other information, such as the covariance of the output and state estimation error in the case of time-varying Kalman filters are used.
As the goal is to use $i^j_k$, together with $u_k$, $y_k$, to estimate $\rho_k$, it is important to 
correctly tune the observers 
associated with the $N_\theta$ models in $\Theta$ to ensure that each $i^j_k$ 
is meaningful. For example, a slower observer may be more robust against 
measurement noise, but its ``inertia'' in reacting to changes may compromise 
the effectiveness of the resulting virtual sensor. 

The computational burden introduced by the observers also needs to be considered. As it will be necessary to run the full bank of $N_\theta$ observers in parallel in real-time, 
a viable option is to use the standard Luenberger observer~\cite{ModernControlSystemShinners1998} 
\begin{equation}
\left\{
\begin{array}{rl}
\hat \xi^j_{k+1}=&  A_j\hat \xi_k^j +d_j + B_j u_k- L_j (\hat y^j_k -y_k)\\
\hat y^j_k=& C_j \hat x_k^j + e_j
\end{array}
\right.
\label{eq:observerDynamic}
\end{equation}
where $L_j$ is the observer gain, and set $i^j_k=\hat \xi^j_k$. Since minimal 
state-space realizations are used to define $\Sigma_j$, each pair $(A_j, C_j)$ 
is fully observable, and the eigenvalues of $A_j-L_jC_j$ can arbitrarily be 
placed inside the unit circle. Note also that any technique for choosing $L_j$ 
can be employed here, such as stationary Kalman filtering.

\subsection{A model-free hypothesis testing algorithm}
\label{sec:datadriven}
After the $N_\theta$ observers have been synthesized, we now build a 
hypothesis testing scheme based on them using a \emph{discriminative} 
approach~\cite{Liu2010GenerativeAndDiscriminative}.
To this end, the initial dataset $\mathcal{D}$ is processed to generate the 
information vectors $i^j_k$, $k\in[k_1,\Nsamp]$. 
Let $D_{\rm aug}\coloneqq \{i_k^1,\ldots,i^N_k,u_k,y_k,\rho_k\}$, 
$k=k_1,\dots,\Nsamp$,
denote the resulting augmented dataset that will be used to train a predictor $f_\theta:
\rr^v\times\ldots\times \rr^v\times \rr^{n_u}\times\rr^{n_y}\to\rr^{S}$ 
such that
\begin{equation}
\hat \rho_k=f_\theta(i_k^1,\ldots,i_{k-\ell}^1,\dots,i^N_k,\dots,i_{k-\ell}^N,u_k,y_k)
\end{equation}
is a good estimate of $\rho_k$, where $\ell\geq 0$ is a window size
to be calibrated.
Consider the minimization of a loss function $L: \rr^{S}\times\rr^{S} \to \rr$ 
that penalizes the distance between the measured value $\rho_k$ and its 
reconstructed value $\hat \rho_k$, namely a solution of 
\begin{equation}
	 \underset{\theta}{\min} \sum_{k=\ell+1}^\Nsamp 
   L(\rho_k,f_\theta(i_k^1,\dots,i^N_{k-\ell},u_k,y_k))
   \label{eq:trainingProblem}
\end{equation}
    
Solving the optimization problem~\eqref{eq:trainingProblem} directly, however, may be excessively complex, as no additional knowledge about the relation between $\rho_k$ and $\{i_k^1,\ldots,i^N_{k-\ell},u_k,y_k\}$ is taken into account. 
In order to model such a relation, one can rewrite $f_\theta$ as the concatenation of two maps $g_\theta$ and $e_\theta$ such that
\begin{equation}
    \begin{array}{rcl}
    \hat\rho_k&=&g_\theta(\II_k)\\
    \II_k&=&e^\mathrm{FE}\theta(i_k^1,i_{k-\ell}^1,\ldots,i^N_k,\ldots,i_{k-\ell}^N,u_k,y_k)
\end{array}
\label{eq:f}
\end{equation}
where $\II_k$ is a \emph{feature vector} constructed by a given
feature extraction (FE) map $e^\mathrm{FE}: \rr^{v}\times\ldots\times\rr^{v} \times 
\rr^{n_u} \times \rr^{n_y} \to \rr^{n_I}$ from $i_k^1$, $i_{k-\ell}^1$, 
$\ldots$, $i^N_k$, $\ldots$, $i_{k-\ell}^N$, $u_k$, and $y_k$, and 
$g_\theta:\rr^{n_I}\to \rr^S$ is the prediction function to learn from 
the dataset $D_{\rm aug}$.

\begin{figure}
    \begin{center}
        \includegraphics[ trim=0cm 0cm 0cm 0cm, clip,
        width=.8\hsize]{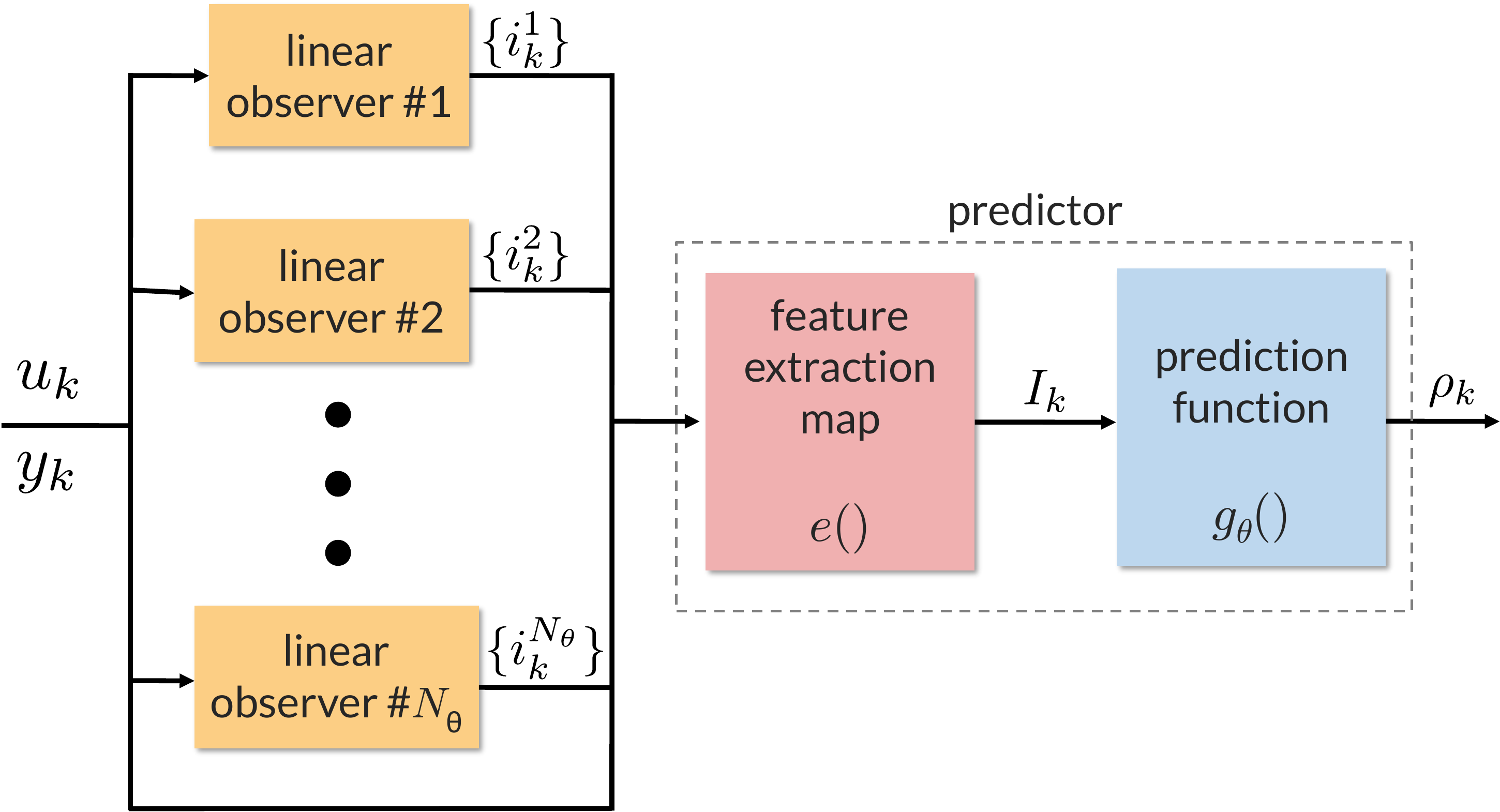}    %
        \caption{Virtual sensor architecture: bank of linear observers, feature extraction map $e$, and prediction function $g_\theta$.} 
        \label{fig:ANNStructure}
    \end{center}
\end{figure}

We propose the following two alternatives for the FE map, namely
\begin{subequations}
\begin{equation}
e^\mathrm{FE}_\theta(\II_k)=\{\hat e_k^1, \hat e^1_{k-\ell}, \ldots, \hat e^N_{k}, \ldots,  
\hat e^N_{k-\ell}, u_k, y_k\}
\label{eq:identifyFeatureExtraction}
\end{equation}
where $\hat e_m^j \eqdef   (\hat y^j_{m} -y_{m})$ and, to further reduce the number of features, the more aggressive and higher compression FE map
\begin{equation}
e^\mathrm{FE}_\theta(\II_k)=\{\nu^1_k,\dots,\nu^N_k,u_k,y_k\}
\label{eq:compressionFeatureExtraction}
\end{equation}
\label{eq:FeatureExtractionMaps}%
\end{subequations}
where
$$
    \nu^i_k=\sum_{r=k-\ell}^k \sqrt{\frac{m(r-\ell)}{\ell}(\hat y^i_r- 
    y_r)'(\hat y^i_r - y_r)}
$$
and $m : \nn \to \rr$ is an appropriate weighting function. 

The rationale for the maps in~\eqref{eq:FeatureExtractionMaps} is that one of 
the most common features used in hypothesis testing algorithms is the estimate 
of the covariance of the residuals produced by each observer. Thus, this approach can 
be thus considered a generalization of the window-based 
hypothesis-testing algorithms explored in the literature, such as 
in~\cite{alsuwaidan2011generalized,hanlon2000multiple}.
The FE map~\eqref{eq:compressionFeatureExtraction} brings this idea one step further so that $e^\mathrm{FE}_\theta(\mathcal{I}_k)$ has only $n_I=(N+1)n_y+n_u$ components. This means that the input to $g_\theta$, and therefore the predictor itself, can get very compact. 

Note that both the window-size $\ell$ and the weighting function $m$ are 
hyper-parameters of the proposed approach. In particular, the value of $\ell$ 
must be chosen 
carefully: if it is too small the time window of past output prediction errors 
may not be long enough for a slow observer. On the other hand, if $\ell$ is too 
large the virtual sensor may become excessively slow in detecting changes of 
$\rho_k$. The weighting function $m$ acts in a similar manner and can be used 
for fine-tuning the behavior of the predictor.

Note also that our choices for $e^\mathrm{FE}_\theta$ in~\eqref{eq:FeatureExtractionMaps} 
are not the only possible ones, nor necessarily the optimal ones.
For example, by setting $e^\mathrm{FE}_\theta(\II_k) = \II_k$, and therefore 
$f_\theta=g_\theta$, one recovers the general case 
in~\eqref{eq:trainingProblem}. Finally, note that our analysis has been 
restricted to a pre-assigned 
function $e^\mathrm{FE}$, although this could also be learned from data. To this end, 
the interested reader is referred 
to~\cite{masti2018Autoencoder,Gao2014ExtractFeaturesviaDAE,guyon2006introduction}
 and the references therein.

\subsubsection{Choice of learning techniques}
\label{ANN}

As highlighted in~\cite{Masti2019LearningBinary}, in order to target an embedded implementation, it is necessary to envision a learning architecture for $g_\theta$ that has a limited memory footprint and requires a small and well predictable throughput. 
To do so, instead of developing an application-specific functional approximation  scheme, we resort to well-understood machine-learning techniques. 
In particular, three possible options are explored in this work , all well suited for our purposes and which require a number of floating-point operations (flops) for their evaluation which is independent of the number of samples used in the training phase, in contrast for example to K-nearest neighbor regression~\cite{hastie2009elements}.

\emph{Remark:} Such choices are not the only possible ones and other approaches may be better 
suitable for specific needs. For example, if one is interested in getting an 
uncertainty measure coupled to the predictions, the use 
of regression techniques based on Gaussian processes~\cite{rasmussen2003gaussian}
could be more suitable.

\subsubsection{Compact artificial neural networks} 
Artificial neural networks (ANN) are a widely used machine-learning technique
that has already shown its effectiveness in other MMAE-based 
schemes~\cite{vincent1998MMAEwithANN}.
An option to make ANN very lightweight is to resort on very compact feed-forward topologies comprised of a 
small number of layers and a computationally cheap activation function in their 
hidden neurons~\cite{karg2020efficient}, such as the Rectified Linear Unit (ReLU)~\cite{nair2010ReLU}
\begin{equation}
f_{\rm ReLU}(x)=\max\{0,x\}
\label{eq:ReLU}
\end{equation}  
As we want to predict real-valued quantities, 
we consider a linear activation function for the output layer of the network.

\subsubsection{Decision-tree and random-forest regression}

DTRs of limited depth are in general extremely cheap to evaluate yet offer a 
good approximation power~\cite{mastiPippia2020learningHybMPC}. Other advantages 
of DTRs are that they can also work effectively with non-normalized data, they 
can be well interpreted~\cite{marchese2017comparison}, and the contribution 
provided by each input feature is easily recognizable.
The main disadvantage of DTRs is instead that they can suffer from high variance. For this reason, in this work, we also explore the use of random-forest regressors (RFRs)~\cite{breiman1984classification}, which try to solve the issue by bagging together multiple trees at the cost of both a more problematic interpretation and higher computational requirements.

\subsection{Hyper-parameters and tuning procedures}
\label{sec:summaryPara}
The overall architecture of the proposed virtual sensor is shown in 
Figure~\ref{fig:ANNStructure}. Its main hyper-parameters are: 
\begin{enumerate}
\item the number $N_\theta$ of local models to learn from experimental data,
related to the number of leaves of the DTR (see 
Section~\ref{sec:DTRforSegmentation});
\item the order $M$ of the local models (see Section~\ref{sec:linear_models});%
\item the window size $\ell$ of the predictor, which sets the number of 
past/current input features provided to the predictor at each time to produce 
the estimate $\hat \rho_k$ (see Section~\ref{sec:datadriven}).
\end{enumerate}
From a practical point of view, tuning $M$ is relatively easy, as one can use 
as the optimal cost reached by solving~\eqref{eq:MLPV_learning} an indirect 
performance indicator to properly trade-off between the quality of fit and 
storage 
constraints. Feature selection approaches such as the one presented 
in~\cite{breschiMejari2020shrinkage} can also be used.
The window size $\ell$ of the predictor is also easily tunable by using any 
feature selection method compatible with the chosen regression 
technique.
A more interesting problem is choosing the correct number of local models $N_\theta$, especially if one considers that MMAE-like schemes often do not perform well if too many models are considered~\cite{li1999multiple}.
Finally, we mention that the proposed method also requires defining the feature 
extraction map and the predictor structure.

As with most black-box approaches, and considering the very mild assumption
we made on on the system $\Sigma_P$ that generates the data, the 
robustness of the virtual sensor with respect to noise and other sources of 
uncertainty can only be assessed \textit{a posteriori}. For this reason, Section~\ref{sec:synthBench} below 
reports a thorough experimental analysis to assess such robustness 
properties.

\section{Numerical results}
\label{sec:results}
In this section, we explore the performance of the proposed virtual sensor 
approach on a series of benchmark problems. All tests were performed on a PC 
equipped with an Intel Core i7 4770k CPU and 16 GB of RAM.
The models introduced in this section were used \emph{only} to generate the 
training datasets and to test the virtual sensor and are \emph{totally unknown} to 
the learning method. We report
such models to facilitate reproducing the numerical results reported
in this section.

 \subsection{Learning setup}
\label{sec:setup}
All the ANNs involved in learning the virtual sensors were developed in Python 
using the \texttt{Keras} framework~\cite{chollet2015keras} and are composed of 
3 layers (2 ReLU layers with an equal number of neurons and an linear output
layer), with overall $60$ hidden neurons.
The ANNs are trained using the AMSgrad optimization 
algorithm~\cite{Reddi2018AMSgrad}. During training, 5\% of the training set is 
reserved to evaluate the stopping criteria.

Both DTR and RFR are trained using \texttt{scikit-learn}~\cite{sklearn_api} 
and, in both cases, the max depth of the trees is capped to 15. RFRs consist of 
10 base classifiers.
The DTR used to extract the set of local models, as shown in 
Section~\ref{sec:DTRforSegmentation}, is instead constrained to have a maximum 
number of leaves equal to the number $N_\theta$ of linear models considered in 
each test.
The loss function $L_{M_{LPV}}$ used is the well known mean absolute 
error~\cite{SpringerEncyclopedia2011MAE}. In all other cases, the standard 
mean-squared error (MSE)~\cite{SpringerEncyclopedia2011MSE} is considered. 

The performance of the overall virtual sensor is assessed on the testing dataset in terms of the following fit ratio ($\FIT$) and normalized root mean-square error (NRMSE)
\begin{subequations}
\begin{equation}
\FIT=\max\left\{0,1-\frac{\|\rho_\mathcal{T} - \hat \rho_\mathcal{T}\|_2}{\|\bar \rho - \rho_\mathcal{T}\|_2}\right\}
\label{eq:FIT}
\end{equation}
\begin{equation}
\NRMSE=\max\left\{0,1-\frac{\|\rho_\mathcal{T} - \hat \rho_\mathcal{T}\|_2}{\sqrt{\mathcal{T}}|\max(\rho_\mathcal{T}) - \min(\rho_\mathcal{T})|}\right\}
\label{eq:NRMSE}
\end{equation}
\label{eq:fit-figures}%
\end{subequations}
computed component-wise, where $\bar \rho$ is the mean value of the test sequence $\rho_\mathcal{T}=\{\rho_i\}_{i=1}^\mathcal{T}$ of the true values and $\hat \rho_\mathcal{T}$ is its estimate.

For each examined test case, we report the mean value and standard deviation of 
the two figures in~\eqref{eq:fit-figures} over ten different runs, each one 
involving different realizations of all the excitation signals $u_k$, $p_k$, 
and measurement noise. 

\subsection{A synthetic benchmark system}
\label{sec:synthBench}
We first explore the performance of the proposed approach and analyze the 
effect of its hyper-parameters on a synthetic 
multi-input single-output benchmark problem. Consider the nonlinear 
time-varying system
\begin{subequations}
\begin{equation}
\Sigma_\mathrm{S}=\left\{
\begin{array}{rcl}
x_{k+1}&=&Hx_k+\frac{\alpha}{2}\atan(x_k)+ \log(\rho_k+1)Fu_k \\
\rho_{k+1}&=&h(\rho_k,u_k,k)		\\
y_k&=& -(1+e^{\rho_k})\matrice{ccccc}{0&0&0&0&1} x_{k}
\end{array}
\right.
\label{eq:synthBenchmarkExample}
\end{equation}
where $x\in \rr^5$, $\atan$ is the arc-tangent element-wise operator, 
$\alpha,\rho_k \in \rr$, matrices $H$ and $F$ are defined as

\begin{equation}
\begin{array}{lr}

H={\begin{bmatrix}

 0.0& 0.1&0.0& 0.0& 0.0& \\
 0.0& 0.0&-1.0& 0.0& 0.0& \\
 0.0& 0.0& 0.0& 1.0& 0.0& \\
 0.0& 0.0& 0.0& 0.0& 1.0& \\
-0.00909&0.0329& 0.29013 &-1.05376&1.69967 
\end{bmatrix}}\vspace*{.5cm}\\

F={ \begin{bmatrix}
-0.71985 & -0.1985  \\
0.57661 &  0.917661\\
1.68733 & -0.68733 \\
-2.14341 &  2.94341 \\
1.      &  1.      \\
\end{bmatrix}}

\end{array}
\label{eq:matrices}
\end{equation}
and function $h$ is defined by
\begin{eqnarray}
h(\rho_k,u_k,k)&=&\left\{
\begin{array}{l}
p_k~~\mbox{if~}~p_k \in [-0.95,0.95]\\
\frac{p_k}{2}~\mbox{otherwise}
\end{array}\right.
\label{eq:nl_benchmark_f}
\\
p_k&=&0.999\rho_k+0.03\omega_k,\quad \omega_k\sim \mathcal{N}(0,1) 
\label{eq:StochasticProcess}
\end{eqnarray}%
\label{eq:NL_benchmark}%
\end{subequations}
mimicking the phenomenon of a slow parameter drift.
Unless otherwise stated, in the following we consider $\alpha=1$.

Training datasets of various sizes (up to $25,000$ samples) and a dataset of 
$5,000$ testing samples are generated by exciting the benchmark 
system~\eqref{eq:NL_benchmark} with a zero-mean white Gaussian noise input 
$u_k$ with unit standard deviation. All signals are then normalized using the 
empirical average and standard deviation computed on the training set and 
superimposed with a zero-mean white Gaussian noise with a standard deviation of 
$0.03$ to simulate measurement noise.

We consider local ARX linear models involving past $M=5$ inputs and outputs
and solve Problem~\eqref{eq:MLPV_learning} via a fully connected feed-forward 
ANN. In the feature extraction process, we set the window $\ell=7$ on which the 
feature extraction process operates. In all tests we also assume 
$m(i)\equiv 1$, $\forall i\in\mathbb{Z}$.
Unless otherwise noted, we consider deadbeat Luenberger observers, i.e., we 
place the observer poles in $z=0$ using the \texttt{Scipy} 
package~\cite{jones2001scipy}.

\subsection{Dependence on the number $\Nsamp$ of samples}

\begin{table}[b]
\center
\begin{tabular}{l|c|c|c|c|c|c}
{ No. of acquired samples $\Nsamp$} & 5000 & 15000   & 25000 & 5000 & 15000   & 
25000  \\
&\multicolumn{3}{c|}{FE map~\eqref{eq:compressionFeatureExtraction} }&\multicolumn{3}{c}{FE map~\eqref{eq:identifyFeatureExtraction} }\\
\hline 
average FIT~\eqref{eq:FIT}	&0.695	&0.769	&0.781	&0.710	&0.799	&0.820	\\
standard deviation	&0.077	&0.031	&0.026	&0.079	&0.026	&0.020	\\
average NRMSE~\eqref{eq:NRMSE}	&0.934	&0.951	&0.953	&0.937	&0.957	&0.961	\\
standard deviation 	&0.017	&0.004	&0.002	&0.017	&0.003	&0.002	\\
\end{tabular} 
\caption{Accuracy of the virtual sensor using datasets of different size $K$.}
\label{table:ANNwrtDSsize}
\end{table} 

We analyze the performance obtained by the synthesized virtual sensor with 
respect to the number $\Nsamp$ of samples acquired for training during the 
experimental phase.
Assessing the scalability of the approach with respect to the size of the dataset is extremely interesting because most machine learning techniques, and in particular neural networks, often require a large number of samples to be effectively trained. 

Table~\ref{table:ANNwrtDSsize} shows the results obtained by training the 
sensor with various dataset sizes when $N_\theta=5$ observers and an 
ANN predictor are used both using the proposed  high-compression FE 
map~\eqref{eq:compressionFeatureExtraction} and using the map 
in~\eqref{eq:identifyFeatureExtraction}. It is apparent that good results can already
be  obtained with 15,000 samples. With smaller datasets, fit performance 
instead remarkably degrades, especially when using the less aggressive FE 
map~\eqref{eq:identifyFeatureExtraction}.

\subsection{Robustness toward measurement noise}
\label{sec:RobVSnoise}
We analyze next the performance of the proposed approach in the presence of 
various levels of measurement noise.
In particular, we test the capabilities of the virtual sensor 
with $N_\theta=5$ deadbeat observers when trained and tested using data 
obtained 
from~\ref{eq:synthBenchmarkExample} and corrupted with a zero-mean additive 
Gaussian noise with different values of standard deviation $\sigma_N$. As in the other tests, noise 
is applied to the signal after normalization.
The training dataset contains $25,000$ samples.

The results, reported in Table~\ref{table:VsNoise}, show a very similar trend 
for all three functional approximation techniques and, in particular, a very 
steep drop in performance when moving from $\sigma_N=0.03$ to $\sigma_N=0.06$.  
This fact suggests that a good signal-to-noise ratio is necessary to achieve 
good performance with the proposed approach. This finding is not surprising, as our 
method is entirely data-driven, it has more difficulties in filtering noise out
compared to model-based methods.

Performance of ANN, RFR, and DTR is similar to what observed in the previous 
tests.

\begin{table}[!h]
\center
\begin{tabular}{l|cc|c|c}
\multirow{2}{*}{Predictor}& &\multicolumn{3}{c}{ Standard 
deviation $\sigma_N$ of additive noise} \\
&& $0.01$&$0.03$ & $0.06$  \\
\hline 
DTR  & &  0.753~(0.024)&0.716~(0.033)&0.667~(0.046)\\
RFR &&0.801~(0.021)&0.771~(0.027)&0.731~(0.039)\\
ANN & & 0.807~(0.020) &0.781~(0.026)&0.740~(0.036)\\
\end{tabular} 
\caption{Average FIT~\eqref{eq:FIT} (standard deviation) for the three proposed 
learning architectures different sensor noise intensity.}
\label{table:VsNoise}
\end{table} 

\subsection{Dependence on the prediction function}
We analyze the difference in performance between the three proposed learning models
for function $g_\theta$ when using $N_\theta=5$ linear models. The 
corresponding results are reported in Table~\ref{table:ANNvsRFvsDTR}, where it 
is apparent that
as soon as enough samples are available, both RFRs and ANNs essentially perform the same, especially when using the more aggressive FE map~\eqref{eq:compressionFeatureExtraction}. For smaller training datasets, the ANN-based predictor performs slightly better, especially in terms of variance. Regression trees show worst performance but they are still able to produce 
acceptable estimates.

Regarding the results obtained using the FE map~\eqref{eq:identifyFeatureExtraction}, ANNs are remarkably more effective than the other two methods. In particular, while RFRs still show acceptable performance, DTRs fail almost completely.

\begin{table}[!h]
\center
\begin{tabular}{l|c|c|c|c}
\multirow{2}{*}{Predictor}&\multirow{2}{*}{ FE map}& \multicolumn{3}{c}{ Number 
$\Nsamp$ of acquired samples} \\
&& 5000& 15000 & 25000  \\
\hline 
DTR  & \multirow{3}{*}{\eqref{eq:compressionFeatureExtraction}} &  0.624~(0.094)&0.698~(0.042)&0.716~(0.033)\\
RFR &&0.685~(0.104)&0.754~(0.034)&0.771~(0.027)\\
ANN & & 0.695~(0.077) &0.769~(0.031)&0.781~(0.026)\\
\hline
DTR  & \multirow{3}{*}{\eqref{eq:identifyFeatureExtraction}} &  0.412~(0.156)&0.589~(0.054)&0.646~(0.029)\\
RFR &&0.596~(0.180)&0.743~(0.038)&0.768~(0.038)\\
ANN & & 0.710~(0.079) &0.799~(0.026)&0.820~(0.020)\\
\end{tabular} 
\caption{Average FIT~\eqref{eq:FIT} (standard deviation) for the three proposed learning architectures
    for different numbers $K$ of samples in the training dataset.}

\label{table:ANNvsRFvsDTR}
\end{table} 
\subsection{Dependence on the observer dynamics}
\label{sec:PPvsKF}

The dynamics of state-estimation errors heavily depend on the location of the 
observer poles set by the Luenberger observer~\eqref{eq:observerDynamic}, such 
as due to the chosen covariance matrices in the case stationary KFs are used 
for observer design. In this section, we analyze the sensitivity of the 
performance achieved by the virtual sensor with respect to the chosen settings 
of the observer. 

Using $N_\theta=5$ models again, the Luenberger observers were tuned to have 
their poles all in the same location $z\in\mathbb{C}$ inside the unit disk and 
vary such a location in different tests. In addition, we also consider 
stationary Kalman filters designed assuming the following model
\begin{equation}
\left \{ \begin{array}{rcl}
\xi^j_{k+1}&=& A_j \xi_k + B_j u_k + d_j +w_k\\
y^j_{k}&=&C_j \xi_k + e_j + v_k
\end{array}
\right.
\label{eq:KFparameterization}
\end{equation}
where $w_k \sim \mathcal{N}(0,I)$ and $v_k \sim \mathcal{N}(0,\lambda I)$ are 
uncorrelated white noise signals of appropriate dimensions and $\lambda \geq 
0$. 

The resulting virtual sensing performance figures are reported in 
Table~\ref{table:PPandKFComparison} for a training dataset of $\Nsamp=25,000$ 
samples and RFR-based prediction. While performance is satisfactory in all 
cases, fast observer poles allow better performance when pole placement is 
used. Nevertheless, in all but the deadbeat case, KFs provide better 
performance regardless of the chosen covariance term $\lambda$.

\begin{table}[h]
\center
\begin{tabular}{l|c|c|c|c|c|c}
& \multicolumn{3}{c|}{pole placement} &\multicolumn{3}{c}{Kalman filter} \\
\hline 
observer settings & $z=0.0$& $z=0.4$ & $z=0.8$ &$\lambda=1$ & $\lambda=10$ &$\lambda=0.1$  \\
\hline 
average FIT~\eqref{eq:FIT}&0.771	&0.698	&0.467	&0.773	&0.762	&0.772\\
standard deviation &0.027	&0.033	&0.050	&0.026	&0.029	&0.026\\
average NRMSE~\eqref{eq:NRMSE}&0.951	&0.935	&0.886	&0.951	&0.949	&0.951\\
standard deviation &0.002	&0.004	&0.006	&0.002	&0.003	&0.002\\
\end{tabular} 
\caption{Average prediction performance with respect to observer settings.}
\label{table:PPandKFComparison}
\end{table}

\subsection{Dependence on the number $N_\theta$ of local models}
\label{sec:Nmodels}

To explore how sensitive the virtual sensor is with respect to the number 
$N_\theta$ of local LTI model/observer pairs employed, we consider the 
performance 
obtained using $N_\theta=2,3,4,5,7$ local models on the 25,000 sample dataset. 
The results obtained using RFR based virtual sensors are reported in 
Table~\ref{table:performanceWRTNModels} and show that performance quickly 
degrades if too few local models are employed. At the same time, one can also 
note that a large number of models is not necessarily more effective. This finding
suggests that the proposed virtual sensor can be easily tuned by increasing 
the number of models until the accuracy reaches a plateau.

\begin{figure}[htp]
\begin{center}
\includegraphics[ trim=0cm 0cm 0cm 1.1cm, clip,
width=\textwidth]{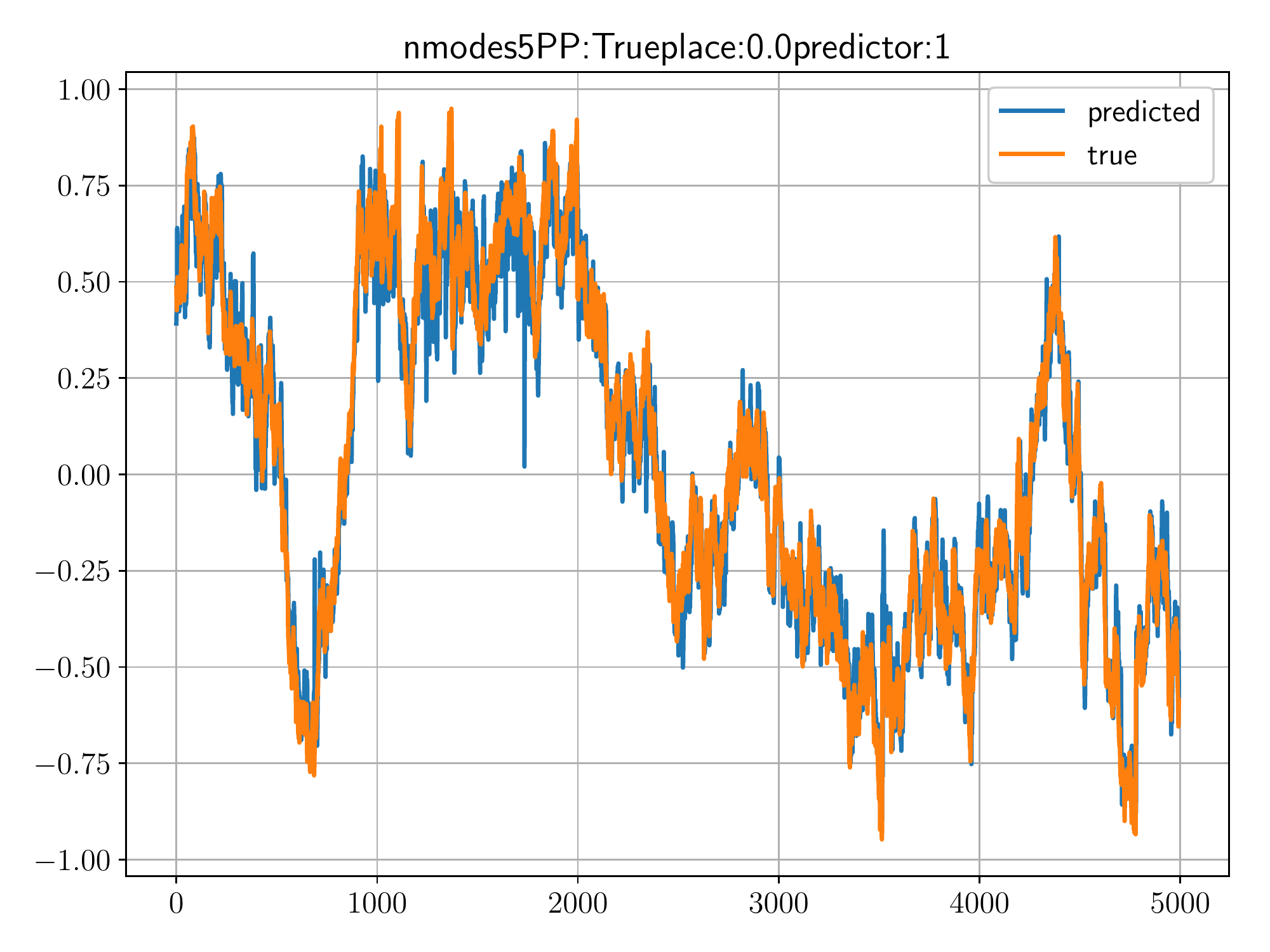}   
\caption{
Example of reconstruction of $\rho_k$ by the virtual sensor based on $N_\theta=5$ local models, using deadbeat observers and a RFR predictor. The figure reports the actual value of $\rho_k$ (orange line) and its estimate $\hat \rho_k$ (blue line).%
}
\label{fig:tracking25k5Example}
\end{center}
\end{figure}

Figure~\ref{fig:tracking25k5Example} shows the estimates of $\rho_k$ obtained by the virtual sensor for $N_\theta=5$, using deadbeat observers and a RFR predictor, for a given realization of~\eqref{eq:StochasticProcess}.

\begin{table}[h]
\center
\begin{tabular}{l|c|c|c|c}
N & 2& 3  & 5  &7 \\ 
\hline 
average FIT~\eqref{eq:FIT}&0.685	&0.766	&0.771	&0.777\\
standard deviation 	&0.035	&0.027	&0.027	&0.028\\
average NRMSE~\eqref{eq:NRMSE}&0.932	&0.950	&0.951	&0.952\\
standard deviation &0.004	&0.002	&0.002	&0.002\\
\end{tabular} 
\caption{Prediction performance of the virtual sensor with respect to the number $N_\theta$ of LTI models.}
\label{table:performanceWRTNModels}
\end{table} 

\subsection{Dependence on the dynamics of  $\rho_k$}
\label{sec:depOnDynamic}
Let us consider a different model 
than~\eqref{eq:nl_benchmark_f}--\eqref{eq:StochasticProcess} to generate the 
value of $p_k$ that defines the signal $\rho_k$, namely the deterministic model
\begin{equation}
p_k=\cos\left(\frac{1}{\beta}k\right)
\label{eq:FMmodulation}
\end{equation}
with $\beta=200$. In this way, it is possible to test the effectiveness of our approach in a purely parameter-varying setting and its robustness against discrepancies between
the way training data and testing data are generated.

The results reported in Tables~\ref{table:resFMLaw} and~\ref{table:resHybrid} are obtained by
using $N_\theta=5$ models and 25,000 samples using all the three different 
prediction functions and deadbeat observers. 
While~\eqref{eq:FMmodulation} is used to generate both training and testing data to produce
the results shown in Table~\ref{table:resFMLaw}, Table~\ref{table:resHybrid} 
shows the case in which the training dataset is generated 
by~\eqref{eq:nl_benchmark_f}--\eqref{eq:StochasticProcess} while the testing 
dataset by~\eqref{eq:FMmodulation}. It is apparent that the proposed approach 
is able to work effectively also in the investigated parameter-varying context 
in all cases and able to cope with sudden changes of $\rho_k$.

\begin{table}[h]
\center
\begin{tabular}{l|c|c|c}
Predictor& DTR & RFR  & ANN  \\
\hline 
average FIT~\eqref{eq:FIT} &0.826	&0.860	&0.859\\
standard deviation 	&0.007	&0.004	&0.005\\
average NRMSE~\eqref{eq:NRMSE}	&0.948	&0.958	&0.957\\
standard deviation	&0.002	&0.001	&0.001\\

\end{tabular} 
\caption{Average accuracy of the virtual sensor employing various kind of prediction functions when both training and testing data are generated by using~\eqref{eq:FMmodulation}.}
\label{table:resFMLaw}
\end{table}

\begin{table}[bh]
\center
\begin{tabular}{l|c|c|c}
predictor& DTR & RFR  & ANN  \\
\hline 
average FIT~\eqref{eq:FIT}&0.755	&0.778	&0.834\\
standard deviation&0.044	&0.066	&0.011\\
average NRMSE~\eqref{eq:NRMSE}		&0.925	&0.932	&0.949\\
standard deviation		&0.014	&0.020	&0.003\\
\end{tabular} 
\caption{Average accuracy of the virtual sensor for different prediction functions
    with training data generated from~\eqref{eq:StochasticProcess} and testing data from~\eqref{eq:FMmodulation}.}
\label{table:resHybrid}
\end{table} 

\subsection{A mode observer for switching linear systems}

\begin{figure}[h]
\begin{center}
\includegraphics[ trim=0cm 0cm 0cm 1.1cm, clip,
width=\textwidth]{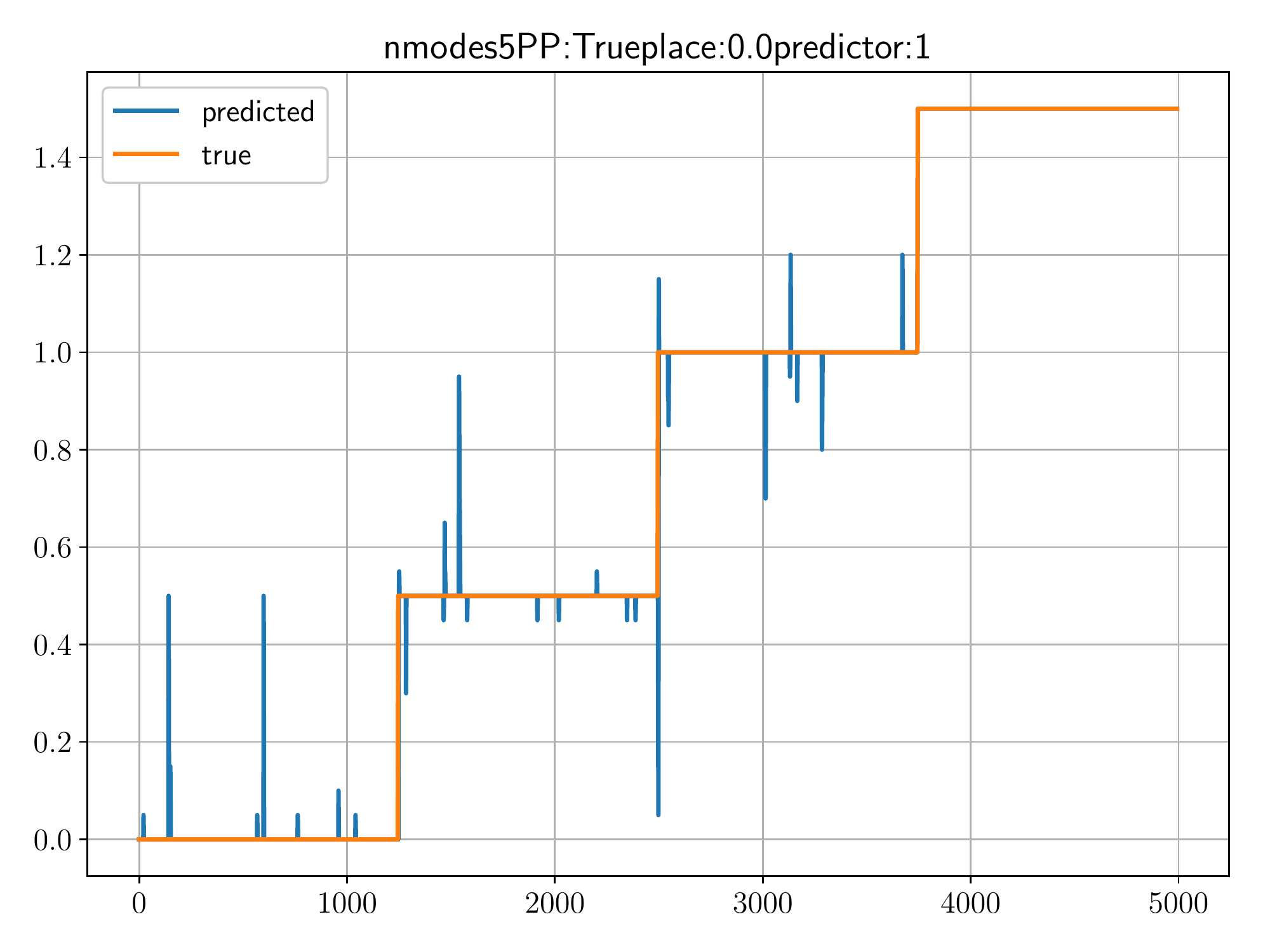}    %
\caption{Mode reconstruction for switching linear systems~\eqref{eq:true-model-SAS}: actual value of the mode $\rho_k$ (orange line) and its estimate $\hat \rho_k$ (blue line) provided by a RFR-based virtual sensor with a bank of 5 deadbeat observers. %
 } 
\label{fig:PWAtrackingRegressor}
\end{center}
\end{figure}

\begin{figure}[h]
\begin{center}
\includegraphics[ trim=0cm 0cm 0cm 1.1cm, clip,
width=\textwidth]{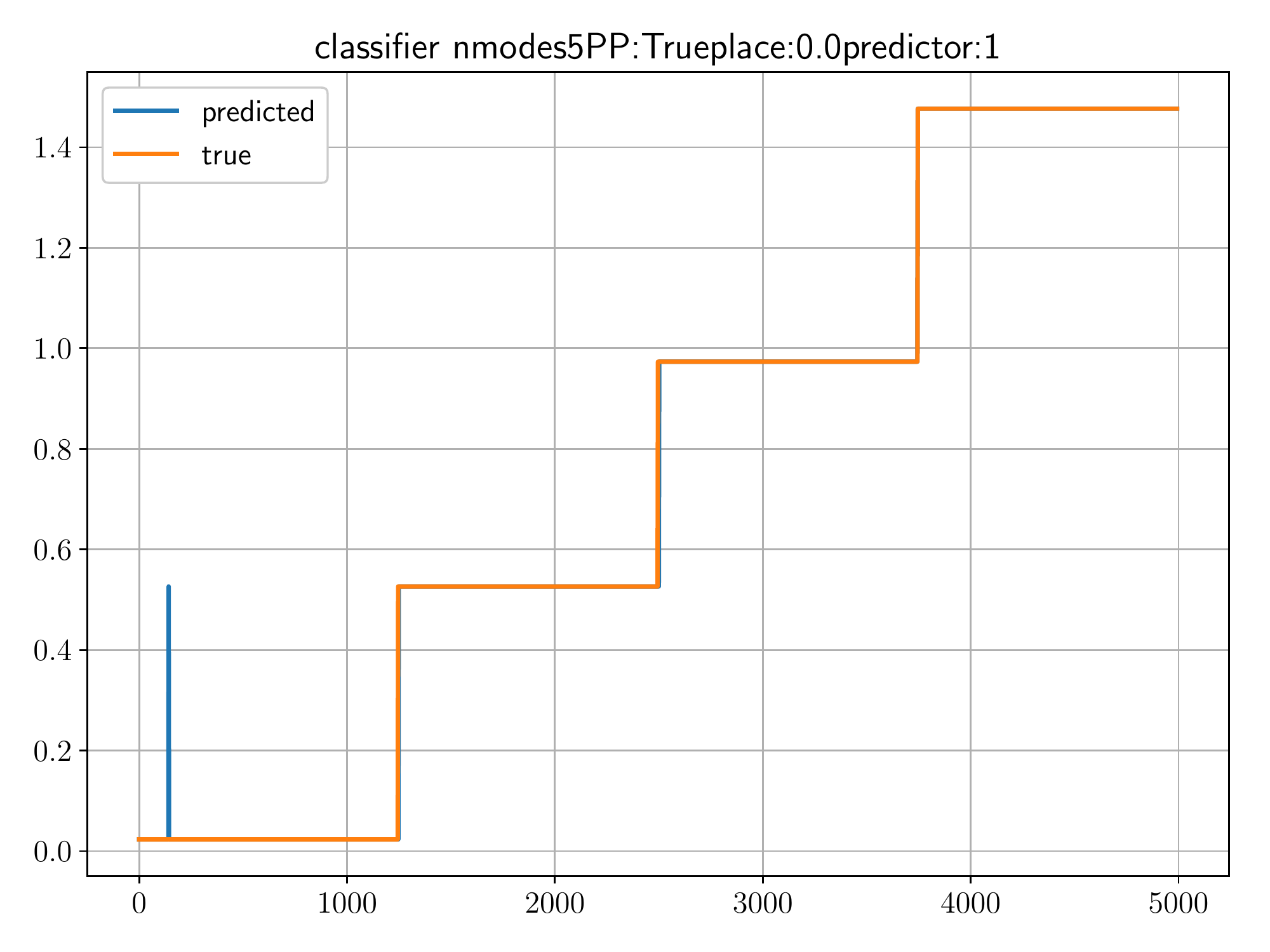}    %
\caption{%
Mode reconstruction for switching linear systems~\eqref{eq:true-model-SAS}: actual value of the mode $\rho_k$ (orange line) and its estimate $\hat \rho_k$ (blue line) provided by a RFC-based virtual sensor and 5 deadbeat observers.
 } 
\label{fig:PWAtrackingClassifier}
\end{center}
\end{figure}

An interesting class of systems which can be described by~\eqref{eq:true-model} are linear switching systems~\cite{liberzon2003switching}, a class of linear parameter varying systems in which  $\rho_k$ can only assume a finite number $s$ of values $\rho^1$, $\ldots$, $\rho^s$. In this case, model~\eqref{eq:true-model} becomes the following discrete-time switching linear system %
\begin{equation}
\Sigma\coloneqq\left\{
\begin{array}{rcl}
x_{k+1}&=&A_{\rho_k}x_k+B_{\rho_k}u_k\\
y_k&=&C_{\rho_k}x_k
\end{array}
\right.
\label{eq:true-model-SAS}
\end{equation}
For switching systems, the problem of estimating $\rho_k$ from input/output measurements is also
known as the \emph{mode-reconstruction} problem.
In order to test our virtual sensor approach for mode reconstruction, we
let the system generating the data be a switching linear system
with $s=4$ modes obtained by the scheduling signal 
\begin{equation}
\rho_{k}=\frac{1}{2}\left\lfloor \frac{4k}{\Nsamp} \right\rfloor
\label{eq:PWAscheduling}
\end{equation}
where  $\Nsamp$ is the number of samples collected in the experiment and 
$\lfloor\cdot\rfloor$ is the downward rounding operator.
For this test, the measurements about the mode acquired during the experiment are noise-free.
As we are focusing on linear switching systems, we set $\alpha=0$.

As in Section~\ref{sec:Nmodels}, we test the performance of the RFR-based virtual sensor trained on 25,000 samples to reconstruct the value of $\rho_k$ when equipped with a different number $N_\theta$ of local models. The corresponding results are reported in Table~\ref{table:PWAresults4} and show that the performance of the sensor again quickly saturates once the number of local models matches the actual number of switching modes.

\begin{table}[h]
    \center
    \begin{tabular}{l|c|c|c|c}
        $N_\theta$& 2 & 3  & 4& 5  \\
        \hline 
        average FIT~\eqref{eq:FIT}&0.785	&0.917	&0.942	&0.943\\
        standard deviation &0.014	&0.016	&0.009	&0.008\\
        average NRMSE~\eqref{eq:NRMSE}&0.920	&0.969	&0.979	&0.979\\
        standard deviation &0.005	&0.006	&0.003	&0.003\\
    \end{tabular} 
    \caption{Accuracy of the virtual sensor employing different predictors for the switching linear system in~\eqref{eq:PWAscheduling}.}
    \label{table:PWAresults4}
\end{table}

The time evolution of the actual mode and the mode reconstructed by the virtual sensor is shown in Figure~\ref{fig:PWAtrackingRegressor}.

\subsubsection{Performance obtained using a classifier in place of a regressor}
The special case of mode reconstruction for switching systems can be also cast 
as a multi-category classification problem. Table~\ref{table:multinomialClass} 
reports the F1-score~\cite{sasaki2007truth} obtained by applying a virtual 
sensor based on a Random Forest Classifier (RFC)  and 5 deadbeat observers to 
discern the current mode of the system. We consider only the case of samples 
correctly labeled, with the RFC subject to the same depth limitation of the 
non-categorical hypothesis tester. 
Table~\ref{table:multinomialClass} also reports the classification accuracy of 
the non-categorical virtual sensor when coupled with a minimum-distance 
classifier (i.e., at each time $k$ the classifier will predict the mode $i$ 
associated with the value $\rho_i$ that is
closest to $\hat \rho_k$). The results refer to a virtual sensor equipped with RFR and 5 deadbeat observers.
It is interesting to note that the classifier architecture is also very 
effective with respect to the FIT metric~\eqref{eq:FIT}, achieving an average 
score of 0.945 with with a standard deviation of  0.010. 

The time evolution of the actual mode and the mode reconstructed by the 
classifier-based virtual sensor is shown in 
Figure~\ref{fig:PWAtrackingClassifier}.

\begin{table}[!h]
\center
\begin{tabular}{r|c|c|c|c}
F1-score / mode \#& 1 & 2  & 3&4  \\
\hline 
RFC&0.997	&0.994	&0.996	&0.998\\ 
standard deviation &0.001	&0.002	&0.002	&0.001\\ 
\hline
RFR&0.996	&0.995	&0.995	&0.997\\
standard deviation &0.002	&0.002	&0.002	&0.002\\ 
\end{tabular} 
\caption{F1-score~\cite{sasaki2007truth} obtained by the RFC-based virtual sensor (RFC) 
    and by the RFR-based virtual sensor + minimum-distance classifier (RFR) 
    on the 4-mode switching linear system~\eqref{eq:PWAscheduling} over 10 runs.}
\label{table:multinomialClass}
\end{table}

\pagebreak

\subsection{Nonlinear state estimation}
\label{sec:Battery}
This sections compares the proposed approach with standard model-based 
nonlinear state-estimation techniques on the problem of 
estimating the state of charge (SoC) of a lithium-ion battery, using the model 
proposed in~\cite{ali2017uas}.

In~\cite{ali2017uas}, the battery is modeled as the following nonlinear third-order dynamical system
\begin{equation}
\Sigma_{\mathrm{Battery}}=\left\{
\begin{array}{l}
\dot x_1(t)= \displaystyle{\frac{\minus i(t)}{C_c}}\\
\dot x_2(t)= \displaystyle{\frac{\minus x_2(t)}{R_{ts}(x_1(t))C_{ts}(x_1(t))} +\frac{i(t)}{C_{ts}(x_1(t))}}\\
\dot x_3(t)= \displaystyle{\frac{\minus x_3(t)}{R_{tl}(x_1(t))C_{tl}(x_1(t))} + \frac{i(t)}{C_{ts}(x_1(t))}}\\
y(t)=E_0(x_1(t))-x_2(t)-x_3(t)-i(t)R_s(x_1(t))
\end{array}
\right.
\label{eq:battery}
\end{equation}
where $x_1(t)$ is the SoC (pure number  $\in [0,1]$ representing the fraction of the battery rated capacity that is available~\cite{Abdi2017SoCDefinitions}), $y(t)$ [V] the voltage at the terminal of the battery, $i(t)$ [A] the current flowing through the battery,  
\begin{equation*}
\begin{array}{l}
 E_0(x_1)=-a_1e^{-a_x x_1}+a_3+a_4 x_1 -a_5 x_1^2 + a_6 x_1^3\\
 R_{ts}(x_1)=a_7e^{-a_8 x_1}+a_9\\
 R_{tl}(x_1)=a_{10}e^{-a_{11} x_1}+a_{12}\\
 C_{ts}(x_1)=a_{13}e^{-a_{17} x_1}+a_{18}\\ R_s(x_1)=a_{19}e^{-a_{20} x_1}+a_{21}
\end{array}
\end{equation*}
and the values of the coefficients $a_{ij}$ correspond to the estimated values reported in Tables~1, 2, 3 of \cite{ali2017uas}.

We analyze the capability of the proposed synthesis method of virtual sensors to reconstruct the value 
$\rho=x_1$ in comparison to a standard extended Kalman filter (EKF)~\cite{anderson1979} based on model~\eqref{eq:battery} and assuming the process noise vector $w_k\in\rr^3$
entering the state equation, $w_{k} \sim \mathcal{N}(0,Q)$, and measurement
noise $v_k \perp w_k \sim \mathcal{N}(0,R)$ on the output $y$ for various realization of $R \in \rr^{3\times 3}, Q \in \rr$. Model~\eqref{eq:battery} is integrated by using an explicit Runge-Kutta 4 scheme.

The simulated system is sampled at the frequency $f_s=\frac{1}{5}$ Hz, starting from a fully charged state $x(0)=[1,~0,~0]'$ and excited with a variable-step current signal $i(t)$ with 
constant amplitude 
\begin{equation}
i(t)=\frac{1}{5}\max\left\{0,\cos\left(\frac{k}{100}\right)+\cos\left(\frac{k}{37}\right)\right\}+u_k
\end{equation}
during the $k$-th sampling steps, with $u_k$ drawn from the uniform distribution
$\mathcal{U}(0,0.4)$.

As the battery will eventually fully discharge, every time the SoC falls below the value $0.05$ the whole state is reset to the initial condition $x(0)$.
The signals $y(t)$ and $i(t)$, once normalized, are processed as described in Section~\ref{sec:setup}. 

For this benchmark, an ANN-based virtual sensor with $N_\theta=5$ local linear models is selected. The corresponding KFs are also designed as described in Section~\ref{sec:PPvsKF} with $\lambda=0.1$, and FE map~\eqref{eq:compressionFeatureExtraction}. Training is performed over 25,000 samples.

\begin{figure}[h]
\begin{center}
\includegraphics[ trim=0cm 0cm 0cm 1.1cm, clip,
width=\textwidth]{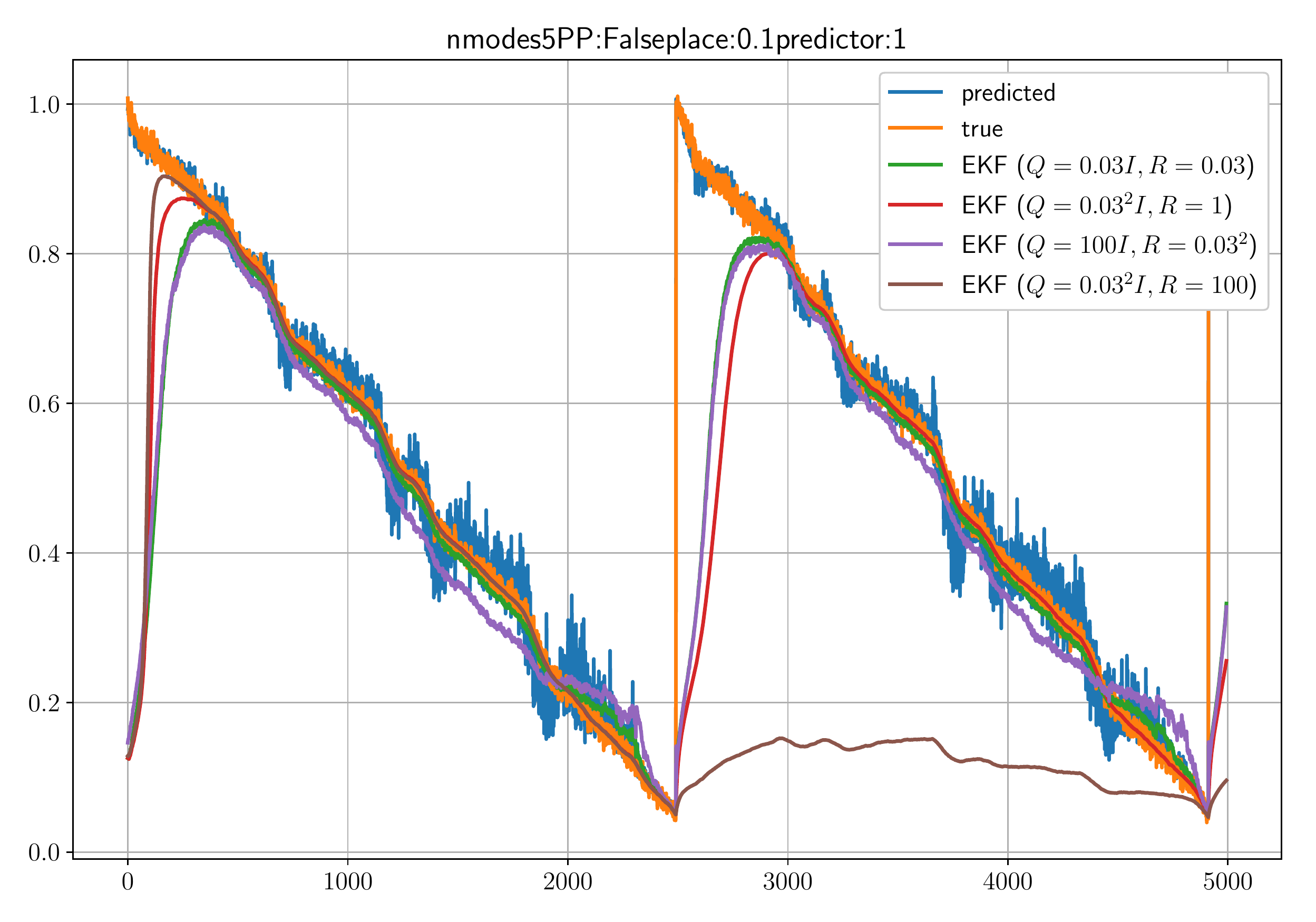}    %
\caption{Estimation of the SoC of the battery: true value $\rho_k$ (orange line), value $\hat \rho_k$
    estimated by the virtual sensor (blue line), values $\hat \rho_k$
    estimated by EKF for different settings of Q and R (green, red, violet, and brown lines). %
     }
\label{fig:EKFvsSensor}
\end{center}
\end{figure}

The results obtained by the virtual sensor and EKFs designed with different
values of the covariance matrices $Q,R$ of process and measurement noise, respectively, are reported in Figure~\ref{fig:EKFvsSensor}. 
While EKF is, in general, more effective in tracking and denoising the true value of the SoC, it performs poorly in terms of bandwidth compared to the proposed virtual sensor, whose performance in terms of filtering noise out remains anyway acceptable.
While both techniques are successful in estimating the SoC of the battery, we remark a main difference between them: EKF \textit{requires} a nonlinear model of the battery, the virtual sensor \textit{does not}.

\subsection{Computation complexity of the prediction functions}

The ANNs used in our tests require approximately between 1,000 and 3,000  weights to be fully parameterized. While this number is fixed by the network topology and depends on the number of inputs to the network, regularization techniques such as $\ell_1$-norm sparsifiers could be used here to reduce the number of nonzero weights (see for instance~\cite{goodfellow2016ANNRegularization} and the reference therein), so to further reduce memory footprint.

Regarding the tree-based approaches, practical storage requirements are 
strongly influenced by the specific implementation and, in general, less 
predictable in advance due to their non-parametric nature. In any case, 
evaluating the predictors on the entire 5,000 sample test set on the reference 
machine only requires a few milliseconds, which makes the approach amenable for 
implementation in most modern embedded platforms.

The training procedure for all the proposed architectures is similarly affordable: on the reference machine, the whole training process is carried out in a few tens of seconds for a training set of 25,000 samples with negligible RAM occupancy.

\section{Conclusions}
\label{sec:conclusions}
This paper has proposed a data-driven virtual sensor synthesis approach, inspired by the MMAE framework, for reconstructing normally unmeasurable quantities such as scheduling parameters in parameter-varying systems, hidden modes in switching systems, and states of nonlinear systems.

The key idea is to use past input and output data (obtained when such quantities were directly measurable) to synthesize a bank of linear observers and use them as a base for feature-extraction maps that greatly simplify the learning process of the hypothesis testing algorithm that estimates said parameters. Thanks to its low memory and CPU requirements, the overall architecture is particularly suitable for embedded and fast-sampling applications.

\bibliographystyle{elsarticle-num}

\subsubsection*{Acknowledgments}

This paper was partially supported by the Italian Ministry of University and Research under the PRIN'17 project ``Data-driven learning of constrained control systems'' , contract no. 2017J89ARP.

\end{document}